\pdfoutput=1 
\documentclass[12pt,reqno]{amsart}

\usepackage{amssymb,amsmath,amsthm,graphicx,xcolor,mathrsfs,tabularx,bbm,tikz,url}
\usepackage{mathabx}\changenotsign
\usepackage{mathtools}
\usepackage{dsfont}
\usepackage[shortlabels]{enumitem}
\usepackage{lmodern}
\usepackage[babel]{microtype}
\usepackage[british]{babel}
\usepackage[utf8]{inputenc}
\allowdisplaybreaks

\usepackage[backref, hypertexnames=false]{hyperref} 
\hypersetup{
	colorlinks,
	linkcolor={red!60!black},
	citecolor={green!60!black},
	urlcolor={blue!60!black}
}
\usepackage{setspace}

\usepackage{geometry}
\let\setminus=\smallsetminus
\let\emptyset=\varnothing

\makeatletter
\def\moverlay{\mathpalette\mov@rlay}
\def\mov@rlay#1#2{\leavevmode\vtop{   \baselineskip\z@skip \lineskiplimit-\maxdimen
		\ialign{\hfil$\m@th#1##$\hfil\cr#2\crcr}}}
\newcommand{\charfusion}[3][\mathord]{
	#1{\ifx#1\mathop\vphantom{#2}\fi
		\mathpalette\mov@rlay{#2\cr#3}
	}
	\ifx#1\mathop\expandafter\displaylimits\fi}
\makeatother

\usepackage{mathfixs}

\usepackage{caption}
\usepackage{subcaption}

\usepackage[mathlines]{lineno}
\usepackage{etoolbox} 

\newcommand*\linenomathpatch[1]{%
	\expandafter\pretocmd\csname #1\endcsname {\linenomath}{}{}%
	\expandafter\pretocmd\csname #1*\endcsname{\linenomath}{}{}%
	\expandafter\apptocmd\csname end#1\endcsname {\endlinenomath}{}{}%
	\expandafter\apptocmd\csname end#1*\endcsname{\endlinenomath}{}{}%
}
\newcommand*\linenomathpatchAMS[1]{%
	\expandafter\pretocmd\csname #1\endcsname {\linenomathAMS}{}{}%
	\expandafter\pretocmd\csname #1*\endcsname{\linenomathAMS}{}{}%
	\expandafter\apptocmd\csname end#1\endcsname {\endlinenomath}{}{}%
	\expandafter\apptocmd\csname end#1*\endcsname{\endlinenomath}{}{}%
}

\expandafter\ifx\linenomath\linenomathWithnumbers
\let\linenomathAMS\linenomathWithnumbers
\patchcmd\linenomathAMS{\advance\postdisplaypenalty\linenopenalty}{}{}{}
\else
\let\linenomathAMS\linenomathNonumbers
\fi

\linenomathpatchAMS{gather}
\linenomathpatchAMS{multline}
\linenomathpatchAMS{align}
\linenomathpatchAMS{alignat}
\linenomathpatchAMS{flalign}
\linenomathpatch{equation}

\usepackage{thmtools} 
\usepackage{cleveref}

\theoremstyle{plain}
\newtheorem{theorem}{Theorem}[section]
\crefname{theorem}{Theorem}{Theorems}

\crefname{proposition}{Proposition}{Propositions}

\crefname{corollary}{Corollary}{Corollaries}

\newtheorem{lemma}[theorem]{Lemma}
\crefname{lemma}{Lemma}{Lemmata}

\newtheorem{conjecture}[theorem]{Conjecture}
\crefname{conjecture}{Conjecture}{Conjectures}

\newtheorem{problem}[theorem]{Problem}
\crefname{problem}{Problem}{Problem}

\crefname{claim}{Claim}{Claims}

\crefname{observation}{Observation}{Observations}

\crefname{setup}{Setup}{Setups}

\crefname{fact}{Fact}{Facts}

\crefname{algorithm}{Algorithm}{Algorithms}

\newtheorem{remark}[theorem]{Remark}
\crefname{remark}{Remark}{Remarks}

\crefname{example}{Example}{Examples}

\theoremstyle{definition}
\newtheorem{definition}[theorem]{Definition}
\crefname{definition}{Definition}{Definitions}

\crefname{construction}{Construction}{Constructions}

\crefname{question}{Question}{Questions}

\numberwithin{equation}{section}

\crefname{section}{Section}{Sections}

\crefname{appendix}{Appendix}{Appendix}

\crefname{figure}{Figure}{Figures}

\def\COMMENT#1{}

\usepackage{comment}

\let\polishlcross=\l
\def\l{\ifmmode\ell\else\polishlcross\fi}

\newcommand{\eps}{\varepsilon}
\renewcommand{\rho}{\varrho}
\newcommand{\sm}{\setminus}
\renewcommand{\subset}{\subseteq}

\newcommand{\NATS}{\mathbb{N}}

\newcommand{\cG}{\mathcal{G}}

\newcommand{\cK}{\mathcal{K}}

\newcommand{\cP}{\mathcal{P}}

\newcommand{\cS}{\mathcal{S}}
\newcommand{\cT}{\mathcal{T}}

\newcommand{\al}{\alpha}
\newcommand{\be}{\beta}
\newcommand{\ce}{\gamma}

\DeclareMathOperator{\hf}{\mathsf{HF}}

\DeclareMathOperator{\hs}{\mathsf{HS}}
\newcommand{\bS}{\mathbb{S}}
\newcommand{\subs}{\subseteq}
\newcommand{\abs}[1]{\left\lvert#1\right\rvert}

\title{Title}
\title[Spanning components and surfaces under minimum vertex degree]{Spanning components and surfaces \\ under minimum vertex degree}

\author[J.~Allsop]{Jack Allsop}

\address[J.~Allsop]{Institut für Mathematik,
	Freie Universität
	Berlin, Germany}
\email{\href{mailto:allsop@mi.fu-berlin.de}{\texttt{allsop@mi.fu-berlin.de}}}

\author[A.~Lamaison]{Ander Lamaison}

\address[A.~Lamaison]{
	Universidad P\'ublica de Navarra, Pamplona, Spain. Previous affiliation: Institute for Basic Science, Daejeon, South Korea
}
\email{\href{mailto:ander.lamaison@unavarra.es}{ander.lamaison@unavarra.es}}

\author[R.~Lang]{Richard Lang}

\address[R.~Lang]{
	Departament de Matemàtiques,
	Universitat Politècnica de Catalunya,
	Barcelona, Spain and
	Centre de Recerca Matemàtica, Barcelona, Spain
}
\email{\href{mailto:richard.lang@upc.edu}{\texttt{richard.lang@upc.edu}}}

\author[S.~Rathke]{Silas Rathke}

\address[S.~Rathke]{
	Institut für Mathematik,
	Freie Universität
	Berlin, Germany
}
\email{\href{mailto:s.rathke@fu-berlin.de}{\texttt{s.rathke@fu-berlin.de}}}

\begin{document}
	
	\begin{abstract}
		We study minimum vertex-degree conditions in $3$-uniform hypergraphs for (tight) spanning components and (combinatorial) surfaces.
		Our main results show that a $3$-uniform hypergraph $G$ on $n$ vertices contains a spanning component if $\delta_1(G) \gtrsim \tfrac{1}{2} \binom{n}{2}$ and a spanning copy of any surface if $\delta_1(G) \gtrsim \tfrac{5}{9} \binom{n}{2}$, which in both cases is asymptotically optimal.
		This extends the work of Georgakopoulos, Haslegrave, Montgomery, and Narayanan who determined the corresponding minimum codegree conditions in this setting.
	\end{abstract}
	

	\maketitle

	\section{Introduction}
	
	Our starting point is Dirac's theorem~\cite{Dir52}, which determines optimal minimum degree conditions that guarantee a Hamilton cycle.
	We study the generalisation of this problem to hypergraphs focusing on the notions of (tight) connectivity and combinatorial surfaces.

	\subsection{Spanning components}\label{ss:SC}
	
	We begin with a discussion of connectivity in hypergraphs.
	To this end, let us introduce some terminology.
	A \emph{$k$-uniform hypergraph} (or \emph{$k$-graph} for short) $G$ consists of a set $V(G)$ of vertices and a set $E(G)$ of edges, where each edge is a $k$-set in $V(G)$.
	We write $e(G) = |E(G)|$.
	For $1 \leq d <k$, the \emph{minimum $d$-degree $\delta_d(G)$} of $G$ is the maximum integer $m$ such that every set of $d$ vertices is in at least $m$ edges. We refer to $\delta_1$ simply as the \emph{minimum (vertex) degree} and to $\delta_{k-1}$ as the \emph{minimum codegree}.
	We remark that minimum $d$-degree conditions are structurally weaker as $d$ decreases.
	It is therefore usually easier to determine optimal codegree conditions than minimum (vertex) degree conditions.
	
	The \emph{line graph} of a $k$-graph $G$ is the $2$-graph on vertex set $E(G)$ with an edge $ef$ whenever $|e \cap f|=k-1$.
	A subgraph of $G$ is \emph{(tightly) connected} if it has no isolated vertices and its edges induce a connected subgraph of its line graph.
	Moreover, we refer to edge-maximal connected subgraphs as \emph{(tight) components}.
	A component is \emph{(vertex) spanning} in $G$ if it contains all vertices of $G$.
	
	It is a classical and simple fact that all $2$-graphs with minimum degree at least $(n-1)/2$ have a spanning component, and that this minimum degree threshold is optimal.
	Similarly, it is a simple exercise to show that a codegree of $n/3$ in a $3$-graph results in a spanning component~\cite{GHM19}, which is also optimal.
	However, despite being a fundamental problem, the minimum codegree threshold for the existence of a spanning component is unknown for $k$-graphs with $k \geq 5$ (see also \cref{con:codegree}).
	In the setting of minimum vertex degrees, it was conjectured by Illingworth, Lang, Müyesser, Parczyk, and Sgueglia~\cite{ILM+24} that a minimum degree of asymptotically $\frac 4 9\binom{n}{2}$ suffices to force a spanning component in a $3$-graph.
	This conjecture turns out to be incorrect as the following lemma shows (see also the construction in \cref{fig:construction}):
	
	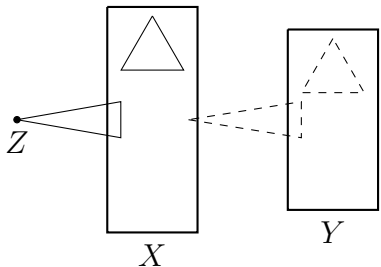
\begin{figure}
		\begin{tikzpicture}[scale=0.6]
			\draw[fill=black, ] (0,0) circle (2pt);
			\draw[thick] (2,2.5) -- (2,-2.5) -- (4,-2.5) -- (4,2.5) -- (2,2.5);
			\draw[thick] (6,2) -- (6,-2) -- (8,-2) -- (8,2) -- (6,2);
			\draw (0,0) -- (2.3,.4) -- (2.3,-.4) -- (0,0);
			\begin{scope}[shift={(3,1.5)}]
				\draw (90:.8) -- (210:.8) -- (330:.8) -- (90:.8);
			\end{scope}
			\begin{scope}[shift={(7,1.)}]
				\draw[dashed] (90:.8) -- (210:.8) -- (330:.8) -- (90:.8);
			\end{scope}
			\draw[dashed] (3.8,0) -- (6.3,-.4) -- (6.3,.4) -- (3.8,0);
			\node at (0,-.5) {$Z$};
			\node at (3,-3) {$X$};
			\node at (7,-2.5) {$Y$};
		\end{tikzpicture}
		\caption{A $3$-graph without a spanning component. The triangles indicate which $3$-sets are edges.
        The edges corresponding to the dashed and non-dashed triangles each form a component.}
		\label{fig:construction}
	\end{figure}
	
	\begin{lemma}\label{lem:HamCompDegOptimal}
		For every sufficiently large integer $n$, there is a $3$-graph $G$ on $n$ vertices with $\delta_1(G) = \frac{1}{2}\binom{n}{2}- O(n)$ that does not have a spanning component.
	\end{lemma}
	
	On the positive side, we are able to determine the (asymptotically) best possible minimum degree threshold for the existence of a spanning component in a $3$-graph.
	The following is our first main result.
	
	\begin{theorem}\label{thm:spanning-component}
		Every $3$-graph $G$ on $n$ vertices with $\delta_1(G) \geq \tfrac{1}{2}\binom{n-1}{2}+1$ contains a spanning component.
	\end{theorem}

	\subsection{Spanning surfaces}
	
	Dirac's famous theorem asserts that any $2$-graph with minimum degree at least $n/2$ contains a spanning cycle.
	Furthermore, this minimum degree threshold is best possible.
	Much work has been done generalising Dirac's theorem to $k$-graphs for $k \geq 3$.
	In the following, we focus on tight Hamiltonicity, which is defined as follows.	
	A \emph{(tight) cycle} $C$ in a $k$-graph $G$ is a subgraph whose vertices can be cyclically ordered such that the edge set of $C$ consists of all sets of $k$ consecutive vertices in this ordering.
	Moreover, $C$ is \emph{Hamilton} if it contains all vertices of $G$.
	
	A natural way to generalise Dirac's theorem to $k$-graphs is to consider minimum codegree thresholds for the existence of a Hamilton cycle.
	This problem has been asymptotically solved for all values of $k$ by R\"odl, Ruci\'nski, and Szemer\'edi~\cite{MinCoDegTHC2}.
	The analogous problem for vertex degrees was solved by Reiher, R\"odl, Ruci\'nski, Schacht, and Szemer\'edi~\cite{RRR19} for $k=3$:
	\begin{theorem}\label{t:mindegHamiltonCycle}
		For every $\eps > 0$, there is a positive integer $n_0$ such that every $3$-graph $G$ on $n \geq n_0$ vertices with $\delta_1(G) \geq (5/9+\eps)\binom{n}{2}$ contains a Hamilton cycle. 
	\end{theorem} 
	More recently, Lang, Schacht, and Volec~\cite{LSV24} determined the vertex degree threshold for $k=4$.
	We note that there are other notions of Hamiltonicity in hypergraphs and results in the same spirit have been proved for them.
	For more details, see, for example, the survey of Simonovits and Szemerédi~\cite{SS19}.
	
	We now turn to a geometric variant of Hamiltonicity.
	A \emph{$(k-1)$-dimensional (combinatorial) sphere} in a $k$-graph $G$ is a subgraph $H \subseteq G$ such that the simplicial complex induced by the edges of $H$ is homeomorphic to the $(k-1)$-dimensional sphere $\mathbb{S}^{k-1}$.
	We say that $H$ is a \emph{(vertex) spanning copy of $\mathbb{S}^{k-1}$} in $G$ if it contains all vertices of $G$. Now assume that $k=3$ and let $\cS$ be a connected, closed $2$-manifold (or \emph{surface} for short). We say that $G$ contains a \emph{spanning copy of $\cS$} if there is a vertex-spanning subgraph $H \subseteq G$ whose edge set induces a simplicial complex that is homeomorphic to~$\cS$.
	
	The edge set of a cycle in a $2$-graph induces a simplicial complex that is homeomorphic to the sphere $\mathbb{S}^1$. This observation motivates another generalisation of Dirac's theorem due to Gowers~(see discussion in \cite{Min2DegSS}): For $1 \leq d < k$, what is the minimum $d$-degree threshold for the existence of a spanning copy of a sphere in a $k$-graph?
	Georgakopoulos, Haslegrave, Montgomery, and Narayanan~\cite{Min2DegSS} recently considered this problem and gave an asymptotically best possible minimum codegree threshold when $k=3$. In fact, they found the minimum codegree threshold for the existence of a spanning copy of a given surface $\cS$ in a $3$-graph.
	
	\begin{theorem}\label{t:min2degSpanningSurface}
		For every surface $\cS$ and $\eps > 0$, there is a positive integer $n_0$ such that any $3$-graph $G$ on $n \geq n_0$ vertices with $\delta_2(G) \geq (1/3+\eps)n$ contains a spanning copy of~$\cS$. 
	\end{theorem}
	
	We remark that a related result in the graph setting was proved by Kühn, Osthus, and Taraz~\cite{KOT05}.
	Recently, Illingworth, Lang, Müyesser, Parczyk, and Sgueglia~\cite{ILM+24} determined minimum supported codegree thresholds for spanning spheres.
    (The minimum supported codegree only ranges over the degrees of sets that are contained in an edge.)
    
	Motivated by \cref{t:min2degSpanningSurface}, we study minimum vertex degree conditions for containing a spanning copy of a given surface.
	Our second main result determines the threshold as follows:
	
	\begin{theorem}\label{thm:spanning-surface}
		For every surface $\cS$ and $\eps > 0$, there exists a positive integer $n_0$ such that every $3$-graph $G$ on $n \geq n_0$ vertices with $\delta_1(G) \geq (5/9 + \eps) \binom{n}{2}$ contains a spanning copy of $\cS$.
	\end{theorem}
	
	Furthermore, the minimum degree condition in \cref{thm:spanning-surface} is asymptotically optimal, as the following lemma shows.
	
	\begin{lemma}\label{lem:5/9-is-optimal}
		For every surface $\cS$ and every $\eps>0$, there is an $n_0$, such that for every $n\geq n_0$ there is a $3$-graph $G$ on $n$ vertices with $\delta_1(G)\ge \frac{5}{9}\binom{n}{2}- \eps n^2$ that does not contain a spanning copy of $\cS$.
	\end{lemma}

	\subsection{Spanning spheres}
	
	Our last main result bounds the thresholds for   spanning spheres by the threshold for Hamilton frameworks.
	Let $\delta_d^{\hs}(k)$ be the infimum over all $\delta \in [0,1]$ such that for every $\eps>0$, there is a positive integer $n_0$ such that every $k$-graph $G$ on $n \geq n_0$ vertices with $\delta_d(G) \geq (\delta + \eps ) \binom{n-d}{k-d}$ contains a spanning copy of the sphere~$\mathbb{S}^{k-1}$.
	We defer the formal definition of the threshold for Hamilton frameworks, denoted by $\delta_d^{\hf}(k)$, to \cref{sec:hamilton-frameworks}.
	For now, we just remark that $\delta_d^{\hf}(k)$  is closely related to the threshold for tight Hamiltonicity, and the two currently agree whenever the latter is known.
	In particular, $\delta_1^{\hf}(3) = 5/9$ in accord with \cref{t:mindegHamiltonCycle}.
	
	\begin{theorem}\label{thm:spanning-spheres-hamilton-framework}
		For every $1 \leq d < k$, we have $\delta_d^{\hs}(k) \leq \delta_d^{\hf}(k)$.
	\end{theorem}

	\subsection*{Organisation of the paper}
	
	This paper is organised as follows. In~\cref{s:HamComp}, we determine the minimum degree threshold for the existence of a spanning component in a $3$-graph by proving \cref{lem:HamCompDegOptimal} and \cref{thm:spanning-component}.
	In~\cref{s:HamSpheres}, we investigate spanning copies of surfaces in $3$-graphs and spanning copies of spheres in $k$-graphs. Our investigations allow us to prove our remaining main results, \cref{thm:spanning-surface}, \cref{lem:5/9-is-optimal}, and \cref{thm:spanning-spheres-hamilton-framework}. Finally, we give some concluding remarks in \cref{s:conc}.

	\section{Spanning components}\label{s:HamComp}
	
	In this section, we prove our main results concerning spanning components in $3$-graphs.
	We begin with the construction for the lower bound.
	
	\begin{proof}[Proof of \cref{lem:HamCompDegOptimal}]
		Let a set $V$ on $n$ vertices be partitioned into three non-empty sets $X$, $Y$, and $Z$.
		Define a $3$-graph $G$ on $V$ by adding all edges of type $XXX$, $XYY$, $YYY$, and $ZXX$ (see \cref{fig:construction}).
		Note that $G$ contains only two components, both of which are not spanning: one using the edges of type $XXX$ and $ZXX$ and another using the edges of type $XYY$ and $YYY$.
        Moreover, if we set $|X| = \lfloor n/\sqrt{2} \rfloor$ and $|Z|=1$, then:
        \begin{itemize}
            \item the degree of a vertex in $X$ is $\binom{|X|-1}{2}+|X|-1+\binom{|Y|}{2} = (2-\sqrt{2})\binom{n}{2}-O(n)$,
            \item the degree of a vertex in $Y$ is $\binom{|Y|-1}{2}+|X|(|Y|-1)=\frac12\binom{n}{2}-O(n)$, and
            \item the degree of a vertex in $Z$ is $\binom{|X|}{2} = \frac12\binom{n}{2}-O(n)$.
        \end{itemize}
        It follows that $G$ has minimum degree $\tfrac{1}{2} \binom{n}{2} - O(n)$.
	\end{proof}
	
	Let $G$ be a $3$-graph.
	The \emph{link graph} $L(x)$ of a vertex $x \in V(G)$ is the $2$-graph on $V(G)$ whose edges are the pairs $yz$ such that $xyz$ is an edge in $G$. Note that $x$ is an isolated vertex in $L(x)$.
	
	\begin{proof}[Proof of \cref{thm:spanning-component}]
		Let $G$ be a $3$-graph on $n$ vertices with $\delta_1(G) \geq \tfrac{1}{2}\binom{n-1}{2}+1$.
		Consider a vertex $x \in V(G)$, and note that its link graph $L(x)$ satisfies
        \begin{equation*}
            \frac{2e(L(x))}{n-1} \geq \frac{n-2}{2} + \frac{2}{n-1} > \frac{n-2}{2}\,.
        \end{equation*}
        This implies that $L(x)$ contains a vertex of degree at least $\tfrac{n-1}{2}$.
		So in particular, $L(x)$ has a unique component of largest order, which we denote by $C_x$.
		Let $\hat C_x \subset G$ contain the edges $xyz$ with $yz \in C_x$.
		Since~$C_x$ is connected in $L(x)$, it follows that $\hat C_x$ is connected in $G$.
		So $\hat C_x$ is contained in a (tight) component of $G$, which we denote by~$C_x^\ast$.
		
		For the sake of contradiction, let us assume that $G$ does not contain a spanning component.
		Let $c$ be the number of distinct components $C_v^\ast$ over all $v \in V(G)$.
		Evidently, $c \geq 2$.
		In the following, we distinguish between the cases $c=2$, $c=3$, and $c\geq 4$.
		
		\medskip
		\textbf{Case 1: $c=2$.}
		Fix a vertex $x$ that does not appear in all $V(C_v^\ast)$. By assumption, $C_x^\ast$ is not a spanning component, so there is a vertex $y\not\in V(C_x^\ast)$. Since $c=2$, every $C_v^\ast$ is either $C_x^\ast$ or $C_y^\ast$. In particular, $x\not\in V(C_y^\ast)$.
		
		Let $V'\coloneqq V(G)\backslash\{x,y\}$, $L'(x)\coloneq L(x)[V']$, and $L'(y)\coloneqq L(y)[V']$. 
		Since $\delta_1(G) \geq \tfrac{1}{2} \binom{n-1}{2} + 1$, the link graph $L(x)$ has at least $\frac{1}{2}\binom{n-1}{2}+1$ edges. Fewer than $\frac{n}{2}$ of them can be incident to $y$ since, otherwise, $y$ would be in the largest component $C_x$ of $L(x)$. Thus, $e(L'(x)) > \frac{1}{2}\binom{n-1}{2}+1-\frac{n}{2}$. The same is true for $L'(y)$. This implies \begin{align*}
			e\big(L'(x)\big)+e\big(L'(y)\big)&>2\cdot\!\left(\frac{1}{2}\binom{n-1}{2}+1-\frac{n}{2}\right)
			=\binom{n-2}{2}\\&=e\big(\overline{L'(x)}\cap\overline{L'(y)}\big)+e\big(L'(x)\big)+e\big(L'(y)\big)-e\big(L'(x)\cap L'(y)\big)\,.
		\end{align*}
		We can conclude 
		\begin{equation}\label{eq: more in cap than the complement}
			e\big(L'(x)\cap L'(y)\big)>e\big(\overline{L'(x)}\cap\overline{L'(y)}\big)\,.
		\end{equation}
		
		We define
		\begin{align*}A &\coloneqq \big\{v\in V'\colon v\in V(C_x), v\in V(C_y)\big\}\\
			B &\coloneqq \big\{v\in V'\colon v\in V(C_x), v\not\in V(C_y)\big\}\\
			D &\coloneqq \big\{v\in V'\colon v\not\in V(C_x), v\not\in V(C_y)\big\}\,.
		\end{align*}
		We denote the sizes of the sets $A$, $B$, and $D$ by $a$, $b$, and $d$, respectively.
		Because of \[b+d=n-2-\abs{V(C_y)}\le \frac{n}{2}-2\le \abs{V(C_x)}-2<\abs{V(C_x)}=a+b,\]
		we have $d< a$.
		
		By \eqref{eq: more in cap than the complement}, there is an edge $\{u,v\}\in E\big(L'(x)\cap L'(y)\big)$, that is,  $\{u,v,x\},\{u,v,y\}\in E(G)$.
		We show now that $u,v\in D$. Suppose $u\in V(C_x)$. Then there is a $w$ such that $\{u,w,x\}\in C_x^\ast$, but then we also have that $\{u,v,x\}$ and $\{u,v,y\}$ are in $C_x^\ast$ which contradicts the assumption that $y\not\in V(C_x^\ast)$. Thus, $u\not\in V(C_x)$ and in the same way one can show that $u\not\in V(C_y)$. Hence, $u\in D$ and, similarly, $v\in D$. In particular, this implies that for any $e\in E\big(L'(x)\cap L'(y)\big)$, its endpoints must both lie in $D$. Thus, $e\big(L'(x)\cap L'(y)\big)\le \binom{d}{2}\le d^2$. 
		
		Next, suppose that $\{u,v\}$ is a pair with $u\in A$ and $v\in D$. We cannot have $\{u,v\}\in L'(x)$, since then $u\in A\subs V(C_x)$ would imply that $v\in V(C_x)$ contradicting $v\in D$. Similarly, we cannot have $\{u,v\}\in L'(y)$. Thus, $\{u,v\}\in E\big(\overline{L'(x)}\cap\overline{L'(y)}\big)$, which implies that
		$e\big(\overline{L'(x)}\cap\overline{L'(y)}\big)\ge ad$. Therefore, \eqref{eq: more in cap than the complement} implies
		\[ad\le e\big(\overline{L'(x)}\cap\overline{L'(y)}\big)<e\big(L'(x)\cap L'(y)\big)\le d^2\,,\]
		that is, $a<d$ which contradicts our previous result $d< a$.
		
		\medskip
		\textbf{Case 2: $c=3$.}
		Fix vertices $x$, $y$, and $z$ such that $C_x^\ast$, $C_y^\ast$, and $C_z^\ast$ are pairwise distinct. Let $C_x'$ be the graph on vertex set $V(G)$ and edge set $E(C_x)\cup\{ux\colon u\in V(C_x)\}$. We define $C_y'$ and $C_z'$ in a similar way. Note that every pair in $\binom{V(G)}{2}$ can be in at most one of $E(C'_x)$, $E(C'_y)$, and $E(C'_z)$ as, otherwise, the components $C_x^\ast$, $C_y^\ast$, and $C_z^\ast$ would not be distinct. In particular $e(C'_x)+e(C'_y)+e(C_z')\le \binom{n}{2}$.
		
		Denote $A \coloneq V(C_x)$, $B\coloneq V(C_y)$, and $C \coloneq V(C_z)$.
		Set $\al \coloneq\frac{\abs A}{n-1}$, $\be\coloneq\frac{\abs B}{n-1}$, and $\ce\coloneq\frac{\abs C}{n-1}$. Note that $\al,\be,\ce\in[\frac12,1]$. In the following, we will often use that for $\delta\in[0,1]$, we have
		\begin{align}\label{eq:pull out}
			\binom{\delta(n-1)}{2}=\delta^2\binom{n-1}{2}-\frac{(1-\delta)\delta}{2}(n-1) \,.
		\end{align}
		
		Now observe that at most $\binom{n-1-\al(n-1)}{2} =\binom{(1-\al)(n-1)}{2}= (1-\alpha)^2\binom{n-1}{2}-\frac{\alpha(1-\alpha)}2(n-1)$ edges of $L(x)$ are not in $C_x$.
		This, together with the fact that $e(L(x))\ge \frac{1}{2}\binom{n-1}{2}+1$, implies $e(C_x) > \big(\tfrac{1}{2} - (1-\al)^2\big)\binom{n-1}{2}+\frac{\alpha(1-\alpha)}2(n-1)$ and, ignoring the second summand,  $e(C'_x)\ge \frac{n-1}{2}+e(C_x)>\frac{n-1}{2}+ \big(\tfrac{1}{2} - (1-\al)^2\big)\binom{n-1}{2}$. Similarly, $e(C'_y)>\frac{n-1}{2}+ \big(\tfrac{1}{2} - (1-\be)^2\big)\binom{n-1}{2}$ and $e(C'_z)>\frac{n-1}{2}+ \big(\tfrac{1}{2} - (1-\ce)^2\big)\binom{n-1}{2}$.
		Therefore,
		\begin{align}
			\binom n2 &\ge e(C_x')+e(C_y')+e(C_z')\notag\\&>\frac{3}{2}(n-1)+\left(\frac{3}{2}-\big(1-\al\big)^2-\big(1-\be\big)^2-\big(1-\ce\big)^2\right)\binom{n-1}{2}\label{equ:betterbound}\\\label{equ:Cx-cup-Cy-cup-Cz}
			&>n-1+\left(\big(2\al-\al^2\big)+\big(2\be-\be^2\big)+\big(2\ce-\ce^2\big)-\frac{3}{2}\right)\binom{n-1}{2}\,.
		\end{align}
		Suppose first that $\al+\be+\ce\geq 2$.
		Then the right side is minimised, for example, when $\al=\tfrac{1}{2}$, $\be=\tfrac{1}{2}$, and $\ce=1$ since $f(a)=2a-a^2$ is both concave and monotone increasing on $\big[\frac{1}{2},1\big]$. However, in this case the right side of \eqref{equ:Cx-cup-Cy-cup-Cz} equals 
		\[n-1+\binom{n-1}{2}=\binom{n}{2}\,,\]
		a contradiction. So we can assume that $\al+\be+\ce< 2$.

        \renewcommand{\bar}[1]{\overline{#1}}
		Let us now find a different upper bound for $e(C'_x)+e(C'_y)+e(C_z')=e(C_x'\cup C_y'\cup C_z')$.
        For this, we define $\bar A=A\cup\{x\}$, $\bar B=B\cup\{y\}$, $\bar C=C\cup\{z\}$. 
		By the principle of inclusion and exclusion, we obtain that
		\begin{align}
			e(C'_x \cup C'_y \cup C'_z) &\leq \binom{|\bar A|}{2} + \binom{|\bar B|}{2} + \binom{|\bar C|}{2} +  \binom{|\bar A \cap \bar B \cap \bar C|}{2}  \nonumber \\ 
			&\quad - \binom{|\bar A \cap \bar B|}{2} - \binom{|\bar B \cap \bar C|}{2} - \binom{|\bar C \cap \bar A|}{2} \,. \label{equ:AcupBcupC}
		\end{align}
		We denote the right side of \eqref{equ:AcupBcupC} by $g(\bar A,\bar B,\bar C)$. We claim that there are sets $A',B',C'$ with $\abs{A'}=\abs{\bar A}$, $\abs{B'}=\abs{\bar B}$, $\abs{C'}=\abs{\bar C}$, where $g(\bar A,\bar B,\bar C)\le g(A',B',C')$ and $A'\cap B'\cap C'=\emptyset$. Indeed, if already $\bar A\cap \bar B\cap \bar C=\emptyset$, then there is nothing to show. Suppose that there is a $v \in \bar A \cap \bar B \cap \bar C$.
		Since $\abs A+\abs B+\abs C=(\al+\be+\ce)(n-1)<2(n-1)$ is an integer, we have $\abs A+\abs B+\abs C\le 2(n-1)-1$.
        So $\abs{\bar A} + \abs{\bar B} + \abs{\bar C} =\abs A+\abs B +\abs C+3\le 2n$.
        By choice of $v$, it follows that there is an element $u$ that is in at most one of the sets $\bar A$, $\bar B$, or $\bar C$.
        Suppose first that $u$ is in one of these sets, say $u \in \bar A$.
		We may then delete~$v$ from $\bar B$ and add $u$ instead.
		It follows that the sizes of $\bar B \cap \bar C$ and $\bar A \cap \bar B \cap \bar C$ each decrease by one, while those of $\bar A$, $\bar B$, $\bar C$, $\bar A \cap \bar B$ and $\bar A \cap \bar C$ stay the same.
		Since $|\bar B \cap\bar  C| \geq |\bar A \cap \bar B \cap \bar C|$, this implies that $g(\bar A, \bar B, \bar C)$ does not decrease when we make this change.
        Now suppose that $u$ is in none of $\bar A$, $\bar B$ and $\bar C$.
        In that case, we also delete~$v$ from $\bar B$ and add $u$ instead.
        So the sizes of $\bar A \cap \bar B$, $\bar B \cap \bar C$ and $\bar A \cap \bar B \cap \bar C$ each decrease by one, while the sizes of $\bar A$, $\bar B$, $\bar C$ and $\bar A \cap \bar C$ stay the same.
        Since $|\bar A \cap \bar B| \geq |\bar A \cap \bar B \cap \bar C|$, the function $g(\bar A, \bar B, \bar C)$ does not decrease under this change either.
        In either case, since $\abs{\bar A\cap \bar B\cap \bar C}$ decreases for each such change, eventually, we end up with the desired sets $A',B',C'$.
		
		Note that $\abs{A' \cap B'} \geq \abs{A'} + \abs{B'} -n=\abs{\bar A}+\abs{\bar B}-n = \abs A+\abs B+2-n=(\al+\be-1)(n-1)+1$ and similar inequalities hold for $\abs{A' \cap C'}$ and $\abs{B' \cap C'}$. Using this, we may bound the right side of \eqref{equ:AcupBcupC} from above by
		\begin{align}
			g(A',B',C') &\le \binom{\al(n-1)+1}{2}+\binom{\be(n-1)+1}{2}+\binom{\ce(n-1)+1}{2}\nonumber
			\\&\quad-\binom{(\al+\be-1)(n-1)+1}{2}\nonumber 
			\\&\quad-\binom{(\be+\ce-1)(n-1)+1}{2} \nonumber
			\\&\quad-\binom{(\ce+\al-1)(n-1)+1}{2}\nonumber\\
            &=\binom{\al(n-1)}{2}+\binom{\be(n-1)}{2}+\binom{\ce(n-1)}{2}+(\al+\be+\ce)(n-1)\nonumber
			\\&\quad-\binom{(\al+\be-1)(n-1)}{2}-(\al+\be-1)(n-1)\nonumber
			\\&\quad-\binom{(\be+\ce-1)(n-1)}{2}-(\be+\ce-1)(n-1)\nonumber
			\\&\quad-\binom{(\ce+\al-1)(n-1)}{2}-(\ce+\al-1)(n-1)\nonumber\\
			&\overset{(\ref{eq:pull out})}= \left(\al^2+\be^2+\ce^2-(\al+\be-1)^2-(\be+\ce-1)^2-(\ce+\al-1)^2\right)\binom{n-1}{2} \nonumber\\
			& \quad+\left(-\frac{(1-\alpha)\alpha}{2}-\frac{(1-\beta)\beta}{2}-\frac{(1-\gamma)\gamma}{2}\right)(n-1)-(\al+\be+\ce-3)(n-1)\nonumber\\
			&\quad +  \frac{(2-\alpha-\beta)(\alpha+\beta-1)}{2} (n-1)\nonumber
			\\&\quad+\frac{(2-\beta-\gamma)(\beta+\gamma-1)}{2} (n-1)\nonumber
			\\&\quad+\frac{(2-\gamma-\alpha)(\gamma+\alpha-1)}{2} (n-1)\nonumber\\
			&=f_1(\alpha,\beta,\gamma)\binom{n-1}{2}+f_2(\alpha,\beta,\gamma)(n-1)\,,\label{e:g_f1_f2}
		\end{align}
		where 
		$$f_1(\alpha,\beta,\gamma)\coloneq\al^2+\be^2+\ce^2-(\al+\be-1)^2-(\be+\ce-1)^2-(\ce+\al-1)^2$$ 
		and 
		\begin{align*}
			f_2(\alpha,\beta,\gamma) & \coloneq -\frac{(1-\alpha)\alpha}{2}-\frac{(1-\beta)\beta}{2}-\frac{(1-\gamma)\gamma}{2}-(\al+\be+\ce-3)
			\\&\quad+\frac{(2-\alpha-\beta)(\alpha+\beta-1)}{2}
			\\&\quad+\frac{(2-\beta-\gamma)(\beta+\gamma-1)}{2}
			\\&\quad+\frac{(2-\gamma-\alpha)(\gamma+\alpha-1)}{2}\,.
		\end{align*}
        Simplifying 
        yields 
        \begin{equation}\label{e:f2_ineq}
            f_2(\al,\be,\ce)=\frac{9}{8}-\frac{1}{8}(2(\al+\be+\ce)-3)^2\le \frac{9}{8} \leq \frac32.
        \end{equation}
        Combining \eqref{equ:betterbound}, \eqref{equ:AcupBcupC}, \eqref{e:g_f1_f2}, and \eqref{e:f2_ineq} gives
		\begin{align*}\label{equ:final}
			\left(\frac{3}{2}-(1-\al)^2-(1-\be)^2-(1-\ce)^2\right)\binom {n-1} 2 < f_1(\alpha,\beta,\gamma)\binom{n-1}{2}\,.
		\end{align*}
        By using the definition of $f_1(\alpha, \beta, \gamma)$, it follows that
        \begin{align*}
        0 &> \left(\frac32-(1-\alpha)^2-(1-\beta)^2-(1-\gamma)^2\right)-\\
        &\hspace{5mm}\left(\alpha^2+\beta^2+\gamma^2-(\alpha+\beta-1)^2-(\beta+\gamma-1)^2-(\alpha+\gamma-1)^2\right) \\
        & = 2(\al-\tfrac{1}{2})(\be-\tfrac{1}{2}) + 2(\be-\tfrac{1}{2})(\ce-\tfrac{1}{2}) + 2(\al-\tfrac{1}{2})(\ce-\tfrac{1}{2})   \,.
		\end{align*}
		However, this contradicts the fact that $\alpha$, $\beta$, and $\gamma$ are each at least $\tfrac{1}{2}$. This contradiction concludes the case $c=3$.
		
		\medskip
		\textbf{Case 3: $c\geq4$.} Following the same steps as in case 2, we can derive the analogue of~\eqref{equ:Cx-cup-Cy-cup-Cz}, with four components instead of three. Thus, there are $\alpha,\beta,\gamma,\eta\in \big[\frac{1}{2},1\big]$ such that
		\[\binom{n}{2}>n-1+\left(\big(2\al-\al^2\big)+\big(2\be-\be^2\big)+\big(2\ce-\ce^2\big)+\big(2\eta-\eta^2\big)-\frac{4}{2}\right)\binom{n-1}{2} \,.\]
		But the right side is minimised for $\alpha=\beta=\gamma=\eta=\frac{1}{2}$, where it is $n-1+\binom{n-1}{2}=\binom{n}{2}$, which is absurd!
	\end{proof}

	\section{Spanning spheres}\label{s:HamSpheres}
	
	In this section, we prove \cref{thm:spanning-surface,thm:spanning-spheres-hamilton-framework,lem:5/9-is-optimal}. 
	For the first two results, we rely on some machinery that allows us to embed blow-ups of cycles (\cref{sec:hamilton-frameworks}) under minimum degree conditions.
	We then continue to show that spanning spheres (and surfaces) can be embedded into suitable blow-ups of cycles (\cref{sec:sphere-into-blow-up,sec:surfaces-proof}).
	We combine these results in the proofs of \cref{thm:spanning-surface,thm:spanning-spheres-hamilton-framework} (\cref{sec:proofs-main-results-surfaces-spheres}).
	We conclude with the construction for \cref{lem:5/9-is-optimal} in \cref{sec:constructions}.

	\subsection{Spanning frameworks}\label{sec:hamilton-frameworks}
    
	The purpose of this subsection is to state a result that allows us to embed vertex-spanning blow-ups of cycles.
	
	A \emph{(tight) path} $P$ in a $k$-graph $G$ is a subgraph whose vertices can be (non-cyclically) ordered such that the edge set of $P$ consists of all sets of $k$ consecutive vertices in this ordering.
	
	Let $G$ be a $k$-graph. A \emph{blow-up} of $G$ is a $k$-graph $G'$ formed by replacing each vertex $x$ of $G$ with a non-empty set $V_x$ and each edge $\{x_1, x_2, \ldots, x_k\}$ with the complete $k$-partite $k$-graph on $V_{x_1}, V_{x_2}, \ldots, V_{x_k}$. We call the vertex sets $V_x$ \emph{clusters} and we say that $x \in V(G)$ has been \emph{blown up by} $|V_x|$. We define the \emph{projection map} $\phi \colon V(G') \to V(G)$ by $\phi(v) = x$ for every $v \in V_x$.
	
	As we shall see, (tight) Hamiltonicity requires three structural features: \emph{connectivity}, \emph{space}, and \emph{aperiodicity}.
	Connectivity has been defined in Section~\ref{ss:SC}. Next, we formalise the conditions of space and aperiodicity.
	
	Let $G$ be a $k$-graph on $n$ vertices.
	A \emph{fractional matching} of $G$ is a function $w\colon E(G) \to [0,1]$ such that $\sum_{e\colon v \in e} w (e) \leq 1$ for every $v \in V(G)$.
	The \emph{size} of $w$ is $\sum_{e \in E(G)} w (e)$.
	We say that $w$ is \emph{perfect} if it has size $n/k$. We call $G$ \emph{spacious} if it has a perfect fractional matching.
	
	A \emph{homomorphism} between $k$-graphs $C$ and $G$ is a function $\phi \colon V(C) \to V(G)$ that maps edges to edges.
    A \emph{closed walk} is a homomorphism $\phi\colon C \to G$ of a $k$-uniform cycle $C$.
    We often identify $\phi$ with its image.
    The \emph{order} of $\phi$ is the order of $C$. 
	We call $G$ \emph{aperiodic} if it contains a closed walk whose order is congruent to $1$ modulo $k$.
	
	It is easy to see that any $k$-uniform cycle is connected and spacious.
	Moreover, in the Dirac setting, we are typically interested in families $\cG$ of Hamiltonian $k$-graphs that contain elements of all orders $n \in \NATS$.
	In this case, $\cG$ must necessarily contain $k$-graphs that satisfy aperiodicity (for instance those of order coprime to $k$). 
	Extremal constructions of such families tend to also satisfy a fourth property, which ensures a form of \emph{consistency} in these features across closely related members. This motivates the following definition, which appears in~\cite{LS24a}.
	
	\begin{definition}[Hamilton framework]\label{def:hamilton-framework}
		A family $\cG$ of $n$-vertex $k$-graphs has a \emph{Hamilton framework} $F$ if, for every $G \in \cG$, there is an {$n$-vertex} subgraph $F(G) \subset G$ such that:
		\begin{enumerate}[(F1)]
			\item \label{itm:hf-connected} $F(G)$ is a component, \hfill(connectivity)
			\item \label{itm:hf-matching} $F(G)$ has a perfect fractional matching, \hfill(space)
			\item \label{itm:hf-odd} $F(G)$ contains a closed walk of order $1 \bmod k$, and \hfill(aperiodicity)
			\item \label{itm:hf-intersecting} {$F(H-x) \cup F(H-y)$ is connected for any $(n+1)$-vertex $k$-graph $H$ and $x,y \in V(H)$ such that $H-x, H-y \in \cG$.\hfill(consistency)}
		\end{enumerate}
	\end{definition}
	
	It is not true that every member of a family of $k$-graphs which has a Hamilton framework is Hamiltonian\footnote{
		Consider the $2$-graph obtained by appending an edge to an odd cycle. This graph is connected, has a perfect matching, contains an odd cycle but is not Hamiltonian.
	}.
	However, if the Hamilton framework properties are guaranteed by a sufficiently large minimum degree, one does indeed obtain Hamiltonicity.
	This is formalised as follows.
	
	\begin{definition}\label{def:ham-fw-threshold}
		Let $\delta_d^{\hf}(k)$ be the infimum over all $\delta \in [0,1]$ such that for every $\eps>0$, there is some $n_0$ such that for every $n \geq n_0$, the family of $n$-vertex $k$-graphs $G$ with $\delta_d(G) \geq (\delta + \eps ) \binom{n-d}{k-d}$ admits a Hamilton framework.
	\end{definition}
	
	The following result was proved by Lang and Sanhueza-Matamala~\cite[Theorem 1.1]{LS24a}.
	
	\begin{theorem}\label{thm:framework-bandwidth-simple}
		For every $1 \leq d < k$ and $\eps >0$, there exist $c$ and~$n_0$ with the following properties.
		Let $G$ be a $k$-graph on $n\geq n_0$ vertices with \[\delta_d(G) \geq \big(\delta_d^{\hf}(k) + \eps \big) \tbinom{n-d}{k-d}.\]
		Let $H$ be an $n$-vertex blow-up of a $k$-uniform cycle with clusters of (not necessarily uniform) size at most $(\log \log n)^c$.
		Then $H \subset G$.
	\end{theorem}
	
	Note that if the clusters of $H$ are all singletons in Theorem~\ref{thm:framework-bandwidth-simple}, then $H$ is a Hamilton cycle in $G$.
	
	The following result~\cite[Theorem 2.8]{LS24a} says that the threshold for Hamilton frameworks is equal to the threshold for Hamiltonicity in all cases where the latter threshold is known~\cite{LS22,LSV24,PRR+20,MinCoDegTHC2}.
	
	\begin{theorem}\label{thm:thresholds}
		We have $\delta_{k-1}^{\hf}(k) = 1/2$, $\delta_{k-2}^{\hf}(k) = 5/9$, and $\delta_{k-3}^{\hf}(k) = 5/8$.
	\end{theorem}
	
	\subsection{Embedding spheres into blow-ups}\label{sec:sphere-into-blow-up}
	
	The goal of this section is to prove Lemma~\ref{lem:spheres-in-blow-ups} below, which was proved by Kühn and Osthus~\cite{KO05} for $k=3$.
	
	\begin{lemma}\label{lem:spheres-in-blow-ups}
		For every $k \geq 3$, there is an $n_0$ such that for every $n\geq n_0$ there is an $n$-vertex blow-up of a $k$-uniform path with clusters of size at most $k+2$ that contains a spanning copy of $\mathbb{S}^{k-1}$.
	\end{lemma}
	
	Our proof relies on the work of Illingworth, Lang, Müyesser, Parczyk, and Sgueglia~\cite{ILM+24}.
	It is known that certain $k$-partite $k$-graphs contain spanning copies of $\bS^{k - 1}$.
	Given positive integers $k$ and $a_1, \dotsc, a_k$, we denote by \emph{$K_k^{(k)}(a_1, \dotsc, a_k)$} the complete $k$-partite $k$-graph with parts of size $a_1, \dotsc, a_k$.
	
	\begin{lemma}[{\cite[Lemma 5.4]{ILM+24}}]\label{lem:partitesphere}
		For any $k \geq 2$, the following $k$-partite $k$-graphs contain spanning copies of $\bS^{k - 1}$\textup{:}
		\begin{enumerate}[label = \textup{(}\alph*\textup{)}]
			\item \label{itm:partitesphere-a} $K_k^{(k)}(2, \dotsc, 2, 2, \ell, \ell)$ for any $\ell \geq 2$,
			\item \label{itm:partitesphere-b} $K_k^{(k)}(2, \dotsc, 2, 3, \ell, \ell)$ for any $\ell \geq 3$.
		\end{enumerate}
	\end{lemma}
	
	We say that a copy $S$ of $\bS^{k - 1}$ in a blow-up of a path $P$ is \emph{doubly edge-covering} if there are two families $\{f_e \colon e \in E(P)\}$, $\{f'_e \colon e \in E(P)\}$ of facets of $S$ such that each family is vertex-disjoint and, for each $e \in E(P)$, $f_e \neq f'_e$ and $\phi(f_e) = \phi(f'_e) = e$. For positive integers $a_1, \dotsc, a_\ell$ (where $\ell \geq k \geq 2$), denote by \emph{$P_\ell^{(k)}(a_1, \dotsc, a_\ell)$} the blow-up of the $k$-uniform path on $\ell$ vertices where the $i$th vertex has been blown-up by $a_i$.
	
	\begin{lemma}[{\cite[Lemma 5.5]{ILM+24}}]\label{lem:thin-path}
		For $k\geq 2$, the blow-up $P_{k + 1}^{(k)}(1, 2, 2, \dotsc, 2, 2, 1)$ contains a spanning copy of $\bS^{k - 1}$ that is doubly edge-covering.
	\end{lemma}
	
	\begin{lemma}[{\cite[Lemma 5.6]{ILM+24}}]\label{lem:growing-path}
		Let $\ell - 1 \geq k \geq 2$.
		If the blow-up  $P_\ell^{(k)}(a_1, \dotsc, a_\ell)$ contains a spanning copy of~$\bS^{k - 1}$ that is doubly edge-covering, then so does $P_{\ell + 1}^{(k)}(1, a_1 + 1, \dotsc, a_{k - 1} + 1, a_k, \dotsc, a_\ell)$.
	\end{lemma}
	
	Recall that the connected sum of two spheres is again a sphere.
	In the language of combinatorial spheres, this translates to the following.
	
	\begin{remark}\label{rmk:gluecommonface}
		Let $\cK$ and $\cK'$ both be combinatorial $(k-1)$-spheres that intersect (exactly) in a common edge $F$. Let $\cK''$ be obtained from the union of $\cK$ and $\cK'$ by deleting $F$.
		From the geometric perspective, this glues together the corresponding spheres along the hole left by $F$, and so $\cK''$ is also a combinatorial $(k-1)$-sphere.
	\end{remark}
	
	We are now ready to prove Lemma~\ref{lem:spheres-in-blow-ups}.
	
	\begin{proof}[Proof of \cref{lem:spheres-in-blow-ups}]
		Let $n$ be sufficiently large with respect to $k$. By \cref{lem:partitesphere}, the $k$-partite $k$-graph $K_k^{(k)}(2, \dotsc, 2, 2, 3, 3)$ on $2k+2$ vertices and the $k$-partite $k$-graph $K_k^{(k)}(2, \dotsc, 2, 3, 3, 3)$ on $2k+3$ vertices each contain a spanning copy of $\bS^{k - 1}$.
		We can thus take vertex-disjoint copies $S_1,\dots,S_{t}$ of these spheres such that
		\[n' = n+kt-\sum_{i=1}^{t}|V(S_i)|\]    
		is divisible by $k$ and $t \leq k$.
		
		Next, consider the blow-up of a path $P_{k + 1}^{(k)}(1, 2, 2, \dotsc, 2, 2, 1)$, which by \cref{lem:thin-path} contains a spanning copy of $\bS^{k - 1}$ that is doubly edge-covering.
		Using \cref{lem:growing-path} iteratively, we can extend this to the blow-up of a path \[G = P_{n'/k+k-1}^{(k)}(1, 2, \ldots, k-1, k, k, \ldots, k, k-1, \ldots, 2, 1)\] on $n'$ vertices, which contains a spanning copy $S$ of $\bS^{k - 1}$.
        Moreover, we take $S$ to be vertex-disjoint from $S_1,\dots,S_t$.
		
		To finish, we use \cref{rmk:gluecommonface} to glue $S_1,\dots,S_t$ to $S$ on facets. Let $\ell = n'/k-(k-1)$ and let $C_1, C_2, \ldots, C_\ell$ be the clusters in $G$ of cardinality $k$. We may assume that $kt \leq \ell$. Fix $i \in \{1, 2, \ldots, t\}$ and let $\mathcal{C}_i = \{C_{(i-1)k+j} \colon 1 \leq j \leq k\}$. Let $e$ be an edge of $S$ such that $e \cap C \neq \emptyset$ for all $C \in \mathcal{C}_i$ and let $e'$ be an edge of $S_i$. Identify $e$ with $e'$ and glue $S_i$ to $S$ on this edge. Let $D_1, D_2, \ldots, D_k$ be the clusters of $S_i$. Then the cluster $C_{(i-1)k+j}$ of~$S$ can be replaced by $C_{(i-1)k+j} \cup D_j$, which increases its size by at most $2$. Altogether, we have added $|V(S_i)|-k$ vertices to $G$. Performing the described operation for all $i \in \{1, 2, \ldots, t\}$ yields a copy of $\mathbb{S}^{k-1}$ which is a spanning subgraph of an $n$-vertex blow-up of a path with clusters of cardinality at most $k+2$.    
	\end{proof}

	\subsection{Embedding surfaces into blow-ups}\label{sec:surfaces-proof}
	
	In this subsection, we improve Lemma~\ref{lem:spheres-in-blow-ups} when $k=3$ by proving that all $2$-dimensional surfaces can be embedded into blow-ups of paths with constant cluster size.
	
	\begin{lemma}\label{lem:surfaces-in-blow-ups}
		For every surface $\cS$, there are integers $m$ and $n_0$ such that for every $n\geq n_0$, there is an $n$-vertex blow-up of a $3$-uniform path $H$ with clusters of size at most~$m$ such that $H$ contains a spanning copy of $\cS$.
	\end{lemma}
	
	\begin{figure}
		\begin{center}
			
			\begin{tikzpicture}[scale = 0.75]
				\foreach \x in {0,2,4,6}
				\foreach \y in {0,2,4}
				\draw [thick] (\x,\y) -- (\x,\y+2);
				
				\foreach \x in {0,2,4}
				\foreach \y in {0,2,4,6}
				\draw [thick] (\x,\y) -- (\x+2,\y);
				
				\foreach \x in {0,2,4}
				\foreach \y in {0,2,4}
				\draw [thick] (\x,\y) -- (\x+2,\y+2);
				
				\node at (0, 6) [inner sep=0.7mm, circle, fill=red!50] {$1$};
				\node at (0, 4) [inner sep=0.7mm, circle, fill=green!50] {$2$};
				\node at (0, 2) [inner sep=0.7mm, circle, fill=blue!50] {$3$};
				\node at (0, 0) [inner sep=0.7mm, circle, fill=red!50] {$1$};
				
				\node at (2, 6) [inner sep=0.7mm, circle, fill=blue!50] {$4$};
				\node at (2, 4) [inner sep=0.7mm, circle, fill=red!50] {$5$};
				\node at (2, 2) [inner sep=0.7mm, circle, fill=green!50] {$6$};
				\node at (2, 0) [inner sep=0.7mm, circle, fill=blue!50] {$4$};
				
				\node at (4, 6) [inner sep=0.7mm, circle, fill=green!50] {$7$};
				\node at (4, 4) [inner sep=0.7mm, circle, fill=blue!50] {$8$};
				\node at (4, 2) [inner sep=0.7mm, circle, fill=red!50] {$9$};
				\node at (4, 0) [inner sep=0.7mm, circle, fill=green!50] {$7$};
				
				\node at (6, 6) [inner sep=0.7mm, circle, fill=red!50] {$1$};
				\node at (6, 4) [inner sep=0.7mm, circle, fill=green!50] {$2$};
				\node at (6, 2) [inner sep=0.7mm, circle, fill=blue!50] {$3$};
				\node at (6, 0) [inner sep=0.7mm, circle, fill=red!50] {$1$};
				
				\draw [thick] (10,0) -- (18,0) -- (18,6) -- (10,6) -- (10,0);
				\draw [thick] (11.3, 3.5) -- (15.2, 1.9) -- (13, 3.9) -- (11.3, 3.5);
				\draw [thick] (10,0) -- (11.3, 3.5) -- (10,6);
				\draw [thick] (10,6) -- (13, 3.9) -- (14,6);
				\draw [thick] (10,0) -- (12, 1.5) -- (13.5, 1.5) -- (15.2, 1.9 ) -- (14,0);
				\draw [thick] (12, 1.5) -- (11.3, 3.5) -- (13.5, 1.5);
				\draw [thick] (12, 1.5) -- (14, 0) -- (13.5, 1.5);
				\draw [thick] (15.2, 1.9 ) -- (18, 6);
				\draw [thick] (14,6) -- (15.3, 4.9) -- (14.9, 3.4) -- (15.2, 1.9 );
				\draw [thick] (15.3, 4.9) -- (13, 3.9) -- (14.9, 3.4);
				\draw [thick] (15.3, 4.9) -- (18, 6) -- (14.9, 3.4);
				\draw [thick] (14,0) -- (16.5, 1.1) -- (17, 2.7) -- (18, 6);
				\draw [thick] (16.5, 1.1) -- (15.2, 1.9 ) -- (17, 2.7);
				\draw [thick] (16.5, 1.1) -- (18, 0) -- (17, 2.7);
				
				\node at (10, 0) [inner sep=0.7mm, circle, fill=red!50] {$1$};
				\node at (14, 0) [inner sep=0.7mm, circle, fill=green!50] {$2$};
				\node at (18, 0) [inner sep=0.7mm, circle, fill=blue!50] {$3$};
				\node at (10, 6) [inner sep=0.7mm, circle, fill=blue!50] {$3$};
				\node at (14, 6) [inner sep=0.7mm, circle, fill=green!50] {$2$};
				\node at (18, 6) [inner sep=0.7mm, circle, fill=red!50] {$1$};
				\node at (11.3, 3.5) [inner sep=0.7mm, circle, fill=green!50] {$4$};
				\node at (15.2, 1.9 ) [inner sep=0.7mm, circle, fill=blue!50] {$5$};
				\node at (13, 3.9 ) [inner sep=0.7mm, circle, fill=red!50] {$6$};
				\node at (12, 1.5) [inner sep=0.7mm, circle, fill=blue!50] {$7$};
				\node at (13.5, 1.5) [inner sep=0.7mm, circle, fill=red!50] {$8$};
				\node at (14.9, 3.4) [inner sep=0.7mm, circle, fill=green!50] {$9$};
				\node at (15.3, 4.9) [inner sep=0.2mm, circle, fill=blue!50] {$10$};
				\node at (17, 2.7) [inner sep=0.2mm, circle, fill=green!50] {$11$};
				\node at (16.5, 1.1) [inner sep=0.2mm, circle, fill=red!50] {$12$};
				
				\node at (3, -1.3) [] {\large $\cT_{9}$};
				\node at (14, -1.3) [] {\large $\cP_{12}$};
			\end{tikzpicture}
		\end{center}
		\label{fig:torus-projective-plane}
		\caption{The edge sets of $\cT_{9}$ and $\cP_{12}$ consist of all the triangles in the respective figures.}\label{fig:degen-col}
	\end{figure}
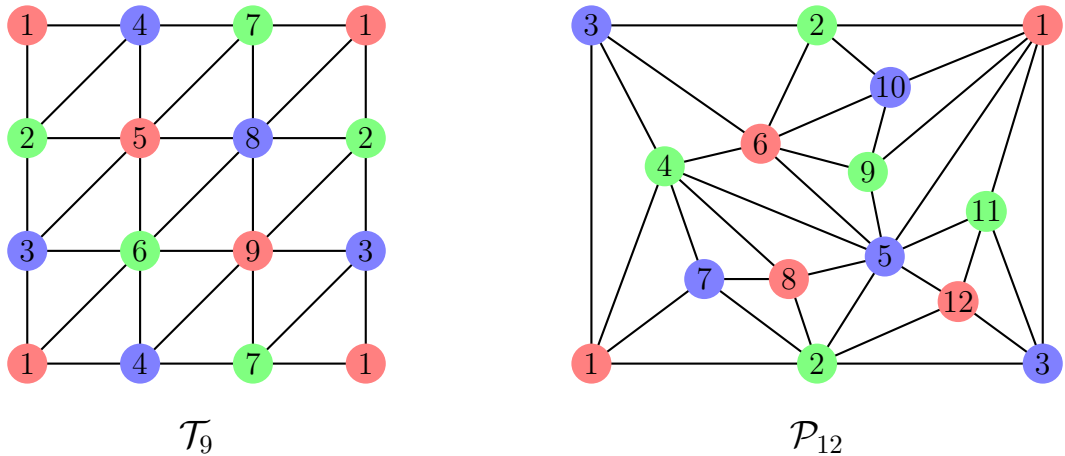
	
	In order to prove Lemma~\ref{lem:surfaces-in-blow-ups}, we will need to use two special $3$-partite $3$-graphs, which we will denote by $\cT_9$ and $\cP_{12}$. These $3$-graphs are defined in Figure~\ref{fig:degen-col}, which appeared in the work of Georgakopoulos, Haslegrave, Montgomery, and Narayanan~\cite{Min2DegSS}. It was noted by the same authors that the simplicial complex induced by $\cT_9$ is homeomorphic to the $2$-dimensional torus and the simplicial complex induced by $\cP_{12}$ is homeomorphic to the real projective plane. We remark moreover that $P_{3}^{(3)}(3,3,3)$ contains a spanning copy of $\cT_9$ and $P_{3}^{(3)}(4,4,4)$ contains a spanning copy of $\cP_{12}$. Given this, and the fact that any $2$-dimensional surface can be obtained from $\mathbb{S}^2$ by gluing some number of two-dimensional tori or real projective planes, the proof of \cref{lem:surfaces-in-blow-ups} follows along the lines of the proof of \cref{lem:spheres-in-blow-ups}.

	\subsection{Embedding spheres and surfaces} \label{sec:proofs-main-results-surfaces-spheres}
	
	It remains to derive \cref{thm:spanning-spheres-hamilton-framework,thm:spanning-surface}.
	
	\begin{proof}[Proof of \cref{thm:spanning-spheres-hamilton-framework}]
        By Dirac's theorem, we may assume that $k\geq 3$.
		Let $\eps>0$ and let $n_0$ be sufficiently large to apply \cref{thm:framework-bandwidth-simple,lem:spheres-in-blow-ups}.
		Let $G$ be a $k$-graph on $n \geq n_0$ vertices with $\delta_d(G) \geq (\delta_d^{\hf}(k) + \eps ) \binom{n-d}{k-d}$.
		Let $H$ be the $k$-uniform sphere on $n$ vertices obtained by \cref{lem:spheres-in-blow-ups}.
		It follows by \cref{thm:framework-bandwidth-simple} that $H \subset G$.
	\end{proof}
	
	Finally, note that \cref{thm:spanning-surface} follows by the combination of \cref{thm:framework-bandwidth-simple,thm:thresholds,lem:surfaces-in-blow-ups}.

	\subsection{Lower bound constructions}\label{sec:constructions}
	
	Finally, we conclude with the proof of \cref{lem:5/9-is-optimal}, which shows that a minimum 1-degree of asymptotically $\frac{5}{9}\binom{n}{2}$ is not sufficient for a spanning copy of a given surface in a $3$-graph.
	
	\begin{proof}[Proof of \Cref{lem:5/9-is-optimal}]
		Consider a partition with non-empty parts $X$ and $Y$ of an $n$-set $V$.
		Let $G$ be the $3$-graph on $V$ obtained by adding all edges of type $XXX$, $XYY$, $YYY$.
		
		Note that $G$ has two components, one consisting of the edges of type $XXX$ and one consisting of the edges of the types $XYY$ and $YYY$.
		We claim that $G$ has no spanning copy of $\cS$ if $|X| > 2|Y|-4 + 4g$, where $g\geq 0$ is the genus of $\cS$.
		Indeed, assume otherwise and let $H \subset G$ be a spanning copy of $\cS$. Since $H$ is vertex-spanning and tightly connected, it must be contained in the component consisting of the edges of the types $XYY$ and $YYY$.
		
		To derive a contradiction, we will use (the generalised) Euler's formula, which states that on the surface $\cS$, the number $f$ of facets, $e$ of 2-edges, and $v$ of vertices satisfy the equation $f+v-e=2-2g$. Since $H$ is a spanning copy of $\cS$, we know that $v=\abs {X}+\abs Y$ and in a triangulation, one always has $e=\frac{3}{2}f$. 
		
		Furthermore, we have $f\ge3\abs{X}$, since every facet has at most one vertex in $X$ and each vertex in $X$ has at least $3$ incident facets. Therefore, we get
		\[2 -2g=f+v-e=-\frac{1}{2}f+v\le -\frac{3}{2}\abs X+\abs X+\abs Y,\]
		which is equivalent to $\abs X\le 2\abs Y-4 +4g$, a contradiction!
		
		Finally, if we take $|X| = \lceil 2n/3\rceil + 4g$, then $G$ satisfies $\delta_1(G) = \tfrac{5}{9} \binom{n}{2}-o(n^2)$.
	\end{proof}

	\section{Conclusion}\label{s:conc}
	
	Connectivity in hypergraphs has thus far mainly been studied as a by-product for more complicated substructures such as Hamilton cycles.
	However, the notion is interesting enough to warrant investigation in its own right.
	We conclude with a few open problems in this direction.
	
	Given Theorem~\ref{thm:spanning-component}, we now know that any $k$-graph $G$ on $n$ vertices with $\delta_1(G) \geq \frac12\binom{n-1}{k-1}+1$ has a spanning component when $k \leq 3$. We believe this remains true for $k > 3$.
	
	\begin{conjecture}
		Let $G$ be a $k$-graph on $n$ vertices with $\delta_1(G) > \tfrac{1}{2}\binom{n-1}{k-1}$.
		Then $G$ contains a spanning component.
	\end{conjecture}
	
	To see that this is best possible, consider the following $k$-graph.
	Let $X$ and $Y$ be a partition of $n-1$ vertices.
	Define a $k$-graph $G'$ on $X \cup Y$ containing all edges except those with $k-1$ vertices in $X$ and one vertex in $Y$.
	Let $G$ be obtained from $G'$ by adding another vertex $z$ and all edges $e$ with $z \in e$ and $e \sm \{z\}\subset X$.
	Then $G$ does not have a spanning component.
	Moreover, by choosing $|X| =  2^{-1/(k-1)}n$, we obtain that $\delta_1(G) \approx \tfrac{1}{2} \binom{n-1}{k-1}$.
	
	We further remark that similar questions can be asked for \emph{$\ell$-components}, where `adjacent' edges intersect in at least $\ell$ vertices.
	(For $\ell=k-1$, this gives the usual tight connectedness.)
	Let us also highlight the following two problems from settings related to (tight) connectivity.
	
	\begin{conjecture}[{\cite[Conjecture 12.4]{LS24a}}]\label{con:tetrahedron}
		Let $G$ be a $3$-graph on $n$ vertices with $\delta_2(G) \geq 3 n / 4$.
		Let $H$ be the $4$-graph on $V(G)$ with an edge for every tetrahedron in $G$.
		Then $H$ is connected.
	\end{conjecture}
	
	A resolution of \cref{con:tetrahedron} would constitute significant progress towards a $3$-uniform version of Pósa's conjecture on squares of cycles~\cite[Section 12.3]{LS24a}.
    A weaker form of the conjecture with $4n/5$ was proved by Bedenknecht and Reiher~\cite{BR20}.
    Very recently, this was improved to $7n/9$ by Bandyopadhyay, Lo and Mycroft~\cite{BLM26}.
	
	\begin{conjecture}[{\cite[Conjecture 5.2]{Min2DegSS}}]\label{con:codegree}
		Every $k$-graph $G$ on $n$ vertices with $\delta_{k-1}(G) > n/k$ contains a spanning component.
	\end{conjecture}
	
	\cref{con:codegree} is trivial for $k=2$ and an easy exercise for $k = 3$~\cite{Min2DegSS}.
    Very recently, the case of $k = 4$ was resolved by Di Braccio, Hearn, Lada, Neve and Zhang~\cite{DHL+26}.
	The cases $k\geq 5$ remain open.
	
	Lastly, a somewhat surprising outcome of our work is that the minimum degree thresholds for Hamilton cycles and surfaces coincide when $k=3$ (see~\cref{t:mindegHamiltonCycle,thm:spanning-surface}).
	Likewise, for $k=2$, this is guaranteed by Dirac's theorem.
	It is a natural question to ask whether this phenomenon extends to higher uniformities.
	By \cref{thm:spanning-spheres-hamilton-framework,thm:thresholds}, we have $\delta_1^{\hs}(4) \leq \delta_1^{\hf}(4) \leq 5/8$  in the $4$-uniform setting.
	However, it is not clear to us that this should be tight.
	
	\begin{problem}
		Prove or disprove that $\delta_1^{\hs}(4) = 5/8$.
	\end{problem}

	\section*{Acknowledgments}
	
	The initial research on this project was conducted during the 2025 Workshop of the Combinatorics and Graph Theory Research Group, FU Berlin at the Sundial Research Center, Zamárdi, Hungary.
	
	\medskip 
	Part of the research on this project was conducted during the visit of the first and the last author at UPC Barcelona supported by the Bilateral AEI+DFG Project PCI2024-155080-2: SRC-ExCo- Structure, Randomness and Computational Methods in Extremal Combinatorics.
	
	\medskip
	Silas Rathke and Jack Allsop were funded by the Deutsche Forschungsgemeinschaft (DFG, German Research Foundation) under Germany's Excellence Strategy – The Berlin Mathematics Research Center MATH+ (EXC-2046/1, project ID: 390685689).
	
	\medskip
	Richard Lang was supported by the Ramón y Cajal programme (RYC2022-038372-I) and by grant PID2023-147202NB-I00 funded
	by MICIU/AEI/10.13039/501100011033.

    \medskip Ander Lamaison was supported by the Institute for Basic Science (IBS-R029-C4).

	\bibliographystyle{plain}
	\bibliography{bibliography.bib}
	
\end{document}